\documentclass[a4paper,11pt]{amsart}

\usepackage{amsfonts,amsmath,mathrsfs,amssymb,amsthm,amscd,latexsym,amstext,amsxtra}

\usepackage[usenames]{color}
\usepackage[all]{xy}
\usepackage{enumerate}
\xyoption{all}
\usepackage{srcltx} %per fer cerca inversa
\usepackage{amsthm, amsmath, amssymb,latexsym}
\usepackage{amsfonts,epsfig,amscd}
\usepackage[]{fontenc}
\usepackage{xy}
\usepackage{enumerate}
\usepackage[]{fontenc}
\usepackage{xy}
\newtheorem{theorem}{Theorem}[section]
\newtheorem{lemma}[theorem]{Lemma}
\newtheorem{corollary}[theorem]{Corollary}
\newtheorem{proposition}[theorem]{Proposition}
\newtheorem{remark}[theorem]{Remark}
\newtheorem{definition}[theorem]{Definition}

\newtheorem{conj}[theorem]{Conjecture}

\newcommand{\nc}{\newcommand}
\nc{\cH}{{\mathcal H}}
\nc{\cA}{{\mathcal A}}
\nc{\cG}{{\mathcal G}}
\nc{\cC}{{\mathcal C}}
\nc{\cD}{{\mathcal D}}
\nc{\cO}{{\mathcal O}}
\nc{\cI}{{\mathcal I}}
\nc{\cB}{{\mathcal B}}
\nc{\cY}{{\mathcal Y}}
\nc{\cK}{{\mathcal K}}
\nc{\cX}{{\mathcal X}}
\nc{\cS}{{\mathcal S}}
\nc{\cE}{{\mathcal E}}
\nc{\cF}{{\mathcal F}}
\nc{\cZ}{{\mathcal Z}}
\nc{\cQ}{{\mathcal Q}}
\nc{\cN}{{\mathcal N}}
\nc{\cP}{{\mathcal P}}
\nc{\cL}{{\mathcal L}}
\nc{\cM}{{\mathcal M}}
\nc{\cT}{{\mathcal T}}
\nc{\cW}{{\mathcal W}}
\nc{\cU}{{\mathcal U}}
\nc{\cJ}{{\mathcal J}}
\nc{\cV}{{\mathcal V}}
\nc{\bH}{{\mathbb H}}
\nc{\bA}{{\mathbb A}}
\nc{\bG}{{\mathbb G}}
\nc{\bC}{{\mathbb C}}
\nc{\bO}{{\mathbb O}}
\nc{\bI}{{\mathbb I}}
\nc{\bB}{{\mathbb B}}
\nc{\bY}{{\mathbb Y}}
\nc{\bK}{{\mathbb K}}
\nc{\bX}{{\mathbb X}}
\nc{\bS}{{\mathbb S}}
\nc{\bE}{{\mathbb E}}
\nc{\bF}{{\mathbb F}}
\nc{\bZ}{{\mathbb Z}}
\nc{\bQ}{{\mathbb Q}}
\nc{\bN}{{\mathbb N}}
\nc{\bP}{{\mathbb P}}
\nc{\bL}{{\mathbb L}}
\nc{\bM}{{\mathbb M}}
\nc{\bT}{{\mathbb T}}
\nc{\bW}{{\mathbb W}}
\nc{\bU}{{\mathbb U}}
\nc{\bD}{{\mathbb D}}
\nc{\bJ}{{\mathbb J}}
\nc{\bV}{{\mathbb V}}
\nc{\bbZ}{{\mathbb Z}}
\nc{\bR}{{\mathbb R}}
\nc{\fr}{{\rightarrow}}
\nc{\co}{{\nabla}}

\nc{\cu}{{\overlineline{\nabla}}}
\title{Torsion points on theta divisors}
\author{Robert Auffarth, Gian Pietro Pirola and Riccardo Salvati Manni}

\address{R. Auffarth \\Departamento de Matem\'aticas, Facultad de
Ciencias, Universidad de Chile, Santiago\\Chile}
\email{rfauffar@mat.puc.cl }

\address{G. Pirola\\Dipartimento di Matematica,
              Universit\`a di Pavia\\Italy}
\email{gianpietro.pirola@unipv.it}

\address{R. Salvati Manni \\Dipartimento di Matematica ``Guido Castelnuovo'', Universit\`a di Roma ``La Sapienza"\\Italy}
\email{salvati@mat.uniroma1.it}

\thanks{Partially supported by Fondecyt Grant 3150171, Prin 2012 ``Moduli Spaces and Lie Theory'' and Inadm Gnsaga}

\begin{document}

\maketitle

\begin{abstract}
Using the irreducibility of a natural irreducible representation of the theta group of an ample line bundle on an abelian variety, we derive a bound for the number of $n$-torsion points that lie on a given theta divisor. We present also two alternate approaches to attacking the case $n=2$.
\end{abstract}

\section{Introduction}

Let $A$ be a complex abelian variety of dimension $g$ and let $\Theta$ be an ample divisor on $A$ that gives a principal polarization $\mathcal{L}:=\mathcal{O}_A(\Theta)$ (i.e. $\dim H^0(A,\mathcal{L})=1$). We will use the notations $(A,\Theta)$ and $(A,\mathcal{L})$ interchangeably. For $n\geq2$, define 
$$\Theta(n):=\#A[n]\cap \Theta,$$
where $A[n]$ is the group of $n$-torsion points on $A$. It is well-known that $\Theta$ does not contain all $n$-torsion points; this follows easily, for example, from the irreducibility of the representation of the theta group of $\mathcal{L}^n$ in $H^0(A,\mathcal{L}^n)$ as we will discuss below.  It is  a classical result, \cite{M3} that the  evaluation at the $n$-torsion points, $n\geq 4$, of Riemann's theta  function
completely determines the abelian  variety embedded in  $\bP^{n^g-1}$. The image is the intersection of all the quadrics containing the  image of the $n$-torsion points. Moreover  the structure of $A[2]\cap\Theta$ tells us if the principally polarized abelian  variety $(A, \Theta)$ is decomposable, \cite{SM} or is the Jacobian of an hyperelliptic curve, \cite{M2}.    Also recently,   in \cite{DFGPS}
it has  been proved  that  $(A, \Theta)$ is decomposable if and only if the image   of the Gauss  map at the smooth points of $\Theta$ in $A[2]\cap \Theta$ is  contained in a  quadric of  $\mathbb{P}^{g-1}$.

 In \cite{PM}, a bound is obtained for the number of 2-torsion points on a theta divisor. Indeed, they show that $\Theta(2)\leq 4^g-2^g$. However, this bound is far from optimal, and in the same paper it is conjectured that the actual bound is $4^g-3^g$ and is achieved if and only if $(A,\mathcal{L})$ is the polarized product of elliptic curves. One could generalize this and conjecture that for $n$-torsion points the bound should be $n^{2g}-(n^2-1)^g$, with equality if and only if $(A,\mathcal{L})$ is the polarized product of elliptic curves. 

Let $\tau\in\mathcal{H}_g$ be a matrix in the Siegel upper-half space, and for $\delta,\epsilon\in\mathbb{R}^g$ and $z\in\mathbb{C}^g$ define the theta function with characteristics
$$\theta\left[\begin{matrix}\delta\\ \epsilon\end{matrix}\right](\tau,z):=\sum_{m\in\mathbb{Z}^g}\mbox{exp}[\pi i(m+\delta)^t\tau (m+\delta)+2\pi i(m+\delta)^t(z+\epsilon)].$$
When $\delta=\epsilon=0$ we obtain \emph{Riemann's theta function} $\theta(\tau,z):=\theta\left[\begin{matrix}0\\ 0\end{matrix}\right](\tau,z)$; the projection of $\{\theta(\tau,\cdot)=0\}$ to $A_\tau:=\mathbb{C}^g/\tau\mathbb{Z}^g+\mathbb{Z}^g$ gives a symmetric theta divisor (i.e. $-\Theta=\Theta$) that we will denote by $\Theta_\tau$. We remark that any complex principally polarized abelian variety is isomorphic to $(A_\tau,\Theta_\tau)$ for some $\tau\in\mathcal{H}_g$. If we put $\mathcal{L}_\tau:=\mathcal{O}_{A_\tau}(\Theta_\tau)$, it is well-known that the set
$$\left\{\theta\left[\begin{matrix}\delta\\ 0\end{matrix}\right](n\tau,nz):\delta\in\frac{1}{n}\mathbb{Z}^g/\mathbb{Z}^g\right\}$$
is a basis for $H^0(A_\tau,\mathcal{L}_\tau^n)$ and the  set 
 $$\left\{\theta\left[\begin{matrix}\delta\\ \epsilon\end{matrix}\right](\tau,nz):\delta, \epsilon\in\frac{1}{n}\mathbb{Z}^g/\mathbb{Z}^g\right\}$$
is a basis for $H^0(A_\tau,\mathcal{L}_\tau^{n^2}).$ . A simple calculation shows that
$$\theta(\tau,z+\tau\delta+\epsilon)=\lambda(z)\theta\left[\begin{matrix}\delta\\ \epsilon\end{matrix}\right](\tau,z)$$
for some nowhere vanishing function $\lambda$, and it immediately follows that if $\Theta=\Theta_\tau$, then $\Theta(n)$ is exactly the number of vanishing \emph{theta constants} $\theta\left[\begin{matrix}\delta\\ \epsilon\end{matrix}\right](\tau,0)$ for $\delta,\epsilon\in\frac{1}{n}\mathbb{Z}^g/\mathbb{Z}^g$. A similar statement holds if $\Theta$ is the pullback of $\Theta_\tau$ by a translation (i.e. $\Theta$ any theta divisor). If $n=2$ and $4\delta^t\epsilon\equiv1\mbox{ (mod }2)$, then the associated theta constant vanishes, and so $\Theta_\tau(2)\geq 2^{g-1}(2^g-1)$. In fact, this is an equality if $A_\tau\in\mathcal{A}_g\backslash\theta^{\text{null}}$, where $\theta^{\text{null}}$ is the divisor consisting of principally polarized abelian varieties such that one of its symmetric theta divisors has a singularity at a point of order 2.\\

The goal of this paper is to give a stronger bound for $\Theta(n)$. Our main theorem gives the following:

\begin{theorem}
Let $(A,\Theta)$ be a complex principally polarized abelian variety. Then
$$\Theta(2)\leq 4^{g}-g2^{g-1}-2^g$$
and for $n\geq3$
$$\Theta(n)\leq n^{2g}-(g+1)n^g.$$
\end{theorem}

We can make this bound better if $(A,\Theta)$ is decomposable.

After proving this theorem, we present   alternative approaches to attacking the number $\Theta(2)$.  One of these points of view will  give  a better bound than the theorem, in  fact we get

\begin{proposition}
Let $(A,\Theta)$ be a principally polarized abelian variety. Then
$$\Theta(2)\leq 4^{g}-\frac{7^g-1}{3^g-1}$$
 
\end{proposition}

We observe  that the methodologies involved are interesting and different from the original approach, and we believe they will be more useful in the future.\\

In  particular in the last approach that  could produce the conjectural  bound, a matrix $M$ appears, induced  by the Weil pairing between the points of order  two in the abelian  variety.  This  matrix appears also in  other fields of mathematics,  in  coding theory as  the matrix associated to  the Macwilliams identity for the weight enumerator of the codes, cf.\cite{CS} page 103,   and in  the  theory of  Borcherds' additive lifting as the matrix associated to a unitary representation of the  integral metaplectic  group on $\bC[\bF_2^{2g}]$,  cf.\cite{FS}\medskip

\noindent\textit{Acknowledgements:} We would like to thank Sam Grushevsky for reading a preliminary version of the paper and pointing out a counterxample to the original proof of the main theorem. We would also like to thank Corrado de Concini for some helpful discussions.

\section{A bound for $\Theta(n)$}

Since $\mathcal{L}$ is a principal polarization, we have that
$$A[n]=\{x\in A:t_x^*\mathcal{L}^n\simeq\mathcal{L}^n\},$$
where $t_x:A\to A$ denotes translation by $x$. Recall that in this case, the \emph{theta group of} $\mathcal{L}^n$ is a certain central extension of $A[n]$ by $\mathbb{G}_m$ which we will denote by $\mathcal{G}_n$:
$$1\to \mathbb{G}_m\to\mathcal{G}_n\to A[n]\to0.$$
Let $\varphi_n:A\to\mathbb{P}H^0(A,\mathcal{L})$ be the morphism associated to the linear system $|\mathcal{L}^n|$. The vector space $H^0(A,\mathcal{L}^n)$ is an irreducible representation for the theta group $\mathcal{G}_n$ where $\mathbb{G}_m$ acts by scalar multiplication (see \cite[Ch. 4]{Kempf} or \cite[Theorem 2, pg. 297]{Mum}), and we therefore obtain a projective representation 
\begin{equation}\label{rep}\nonumber\rho:A[n]\to\mbox{PGL}(H^0(A,\mathcal{L}^n)).\end{equation}
Because of the irreducibility of the representation, we notice that there is no proper linear subspace of $\mathbb{P}H^0(A,\mathcal{L}^n)$ that is invariant under the action of $A[n]$. Moreover, we have that
$$\rho(x)\cdot\varphi_n(y)=\varphi_n(x+y)$$
for every $x\in A[n]$ and $y\in A$.

Let $H\subseteq A[n]$ be a maximal isotropic subgroup with respect to the commutator pairing associated to the theta group of $\mathcal{L}^n$. We say  that  $H$ is \emph{c-isotropic} if it has  a complementary isotropic subgroup $K$. We remark that any maximal isotropic subgroup is isomorphic to $(\mathbb{Z}/n\mathbb{Z})^g$, as is its complementary isotropic subgroup if it exists. Let $H$ be c-isotropic, let $p:A\to A/H=:A_H$ be the natural projection, and let $q:A_H\to A$ be the inverse isogeny. We have a commutative diagram\\

\centerline{\xymatrix{A\ar[r]^{p}\ar[dr]_{n_A}&A_H\ar[d]^{q}\\&A}}

\vspace{0.3cm}

\noindent where $n_A$ denotes multiplication by $n$ on $A$. By descent theory for abelian varieties, we have that there exists a principal polarization $\mathcal{M}$ on $A_H$ such that $\mathcal{L}^n\simeq p^*\mathcal{M}$ and $q^*\mathcal{L}\equiv \mathcal{M}^n$. We see in this case that $\ker q$ is a maximal c-isotropic subgroup of $A_H[n]$. Let $N$  be a complementary isotropic subspace of $\ker q$.

 Define $\Sigma=q^{-1}(\Theta)\in|q^*\mathcal{L}|$ and for $a\in A_H$, define $\Sigma_a:=\Sigma+a$. For every $b\in A_H[n]$, fix a section $s_b\in H^0(A_H,q^*\mathcal{L})$ such that $\Sigma_b=\mbox{div}(s_b)$.

\begin{lemma}\label{basis}
The set $\{s_b:b\in N\}$ is a basis for $H^0(A_H,q^*\mathcal{L})$.
\end{lemma}
\begin{proof}
We see that for all $a\in \ker q$ and $b\in N$, 
$$\Sigma_{a+b}=\Sigma_b+a=q^{-1}(\Theta+q(b))=\Sigma_b.$$
This means that for all $a\in\ker q$ and $b\in N$, there exists $\lambda_a\in \mathbb{G}_m$ such that $t_a^*s_b=\lambda_a s_b$. In other words, $A_H[n]$ acts on the projective span of $\{s_b:b\in N\}$ in $\mathbb{P}H^0(A_H,q^*\mathcal{L})$. Since the theta group representation is irreducible, we must have that the above set generates the whole space. Moreover, $\#N=\dim H^0(A_H,q^*\mathcal{L})$, and the result follows.
\end{proof}

Let $\varphi_H:A_H\to\mathbb{P}H^0(A_H,q^*\mathcal{L})$ be the morphism associated with $|q^*\mathcal{L}|$.

\begin{definition}
For $H$ a maximal c-isotropic subgroup of $A[n]$, let $c_1+H,\ldots, c_{n^g}+H$ be its cosets (we will assume $c_1=0$). We define the integers
\begin{eqnarray}\nonumber Q_{H,c_i}&:=&\dim_\mathbb{C}\mbox{span}\{\varphi_H(q^{-1}(c_i))\}\\
\nonumber Q_H&:=&\sum_{i=1}^{n^g}Q_{H,c_i}\\
\nonumber Q(n)&:=&\max\{Q_H:H\subseteq A[n]\mbox{ max. c-isotropic subgroup}\}.
\end{eqnarray}

\end{definition}

We can use these numbers to obtain a bound on the number of $n$-torsion points lying on $\Theta$.

\begin{proposition}\label{prop1} Let $(A,\Theta)$ be a principally polarized abelian variety and let $n\geq 2$. Then $\Theta(n)\leq n^{2g}-n^g-Q(n)$.
\end{proposition}
\begin{proof}
We will prove that $\Theta(n)\leq n^{2g}-n^g-Q_H$ for every maximal c-isotropic subgroup $H\subseteq A[n]$. Let $S\subseteq H+c_i$ be a subset with $r\leq Q_{H,c_i}$ elements. We will first prove that $\Theta$ does not contain $(H+c_i)\backslash S$. We see that 
\begin{eqnarray}
\nonumber (H+c_i)\backslash S\subseteq \Theta &\Leftrightarrow& q^{-1}((H+c_i)\backslash S)\subseteq\Sigma\\
\nonumber&\Leftrightarrow& (A_H[n]+d_i)\backslash (\ker q+t_1\sqcup\cdots\sqcup \ker q+t_r)\subseteq\Sigma
\end{eqnarray}
where $q(d_i)=c_i$ and the $t_j$ are chosen so that $S=\{q(t_j):j=1,\ldots,r\}$. Assume this occurs. Now for all $b\in N$, 
$$(A_H[n]+d_i)\backslash(\ker q+t_1+b\sqcup\cdots\sqcup \ker q+t_r+b)\subseteq\Sigma_b.$$
It follows that $q^{-1}(c_i)=\ker q+d_i\subseteq\Sigma_b$ for every $b\notin(\ker q+d_i-t_j)\cap N$. We see then there are $n^{g}-r$ options for $b$. Using Lemma \ref{basis}, this implies that $\varphi_H(q^{-1}(c_i))$ is contained in a linear subspace of $\mathbb{P}H^0(A_H,q^*\mathcal{L})$ of dimension $r-1$, a contradiction with the choice of $r$. Therefore in each coset $c_i+H$, there are at most $n^g-Q_{H,c_i}-1$ points that lie on $\Theta$. By adding everything up we get the bound we were looking for.

\end{proof}

\begin{remark} The proof of the previous proposition is valid over any algebraically closed field of characteristic prime to $n$ and for any theta divisor (i.e. not necessarily symmetric). Moreover, the proposition already gives us a better bound than the one in \cite{PM}. Indeed, there can be at most one $Q_{H,c_{i}}$ equal to 0 (this happens when $(A_H,\mathcal{M})$ is the polarized product of elliptic curves), and so $\Theta(2)\leq 4^g-2^g-(2^g-1)=4^g-2^{g+1}+1$.
\end{remark}

The next proposition shows that when looking for a bound for $\Theta(n)$, we can always assume that $\Theta$ is given by the zero set of Riemann's theta function.

\begin{proposition}
If $\Theta$ is Riemann's theta function, then equality holds in Proposition \ref{prop1}. 
\end{proposition}
\begin{proof}
Assume that $\Theta$ is Riemann's theta function, and so $\Theta=\Theta_\tau$ on $A_\tau$ for some $\tau\in\mathcal{H}_g$. Let $\Lambda_\tau$ be the lattice $\tau\mathbb{Z}^g+\mathbb{Z}^g$ and take the maximal c-isotropic subgroup $H=\{\tau\epsilon:\epsilon\in\frac{1}{n}\mathbb{Z}^g/\mathbb{Z}^g\}+\Lambda_\tau$ of $A_\tau[n]$. We have the quotient maps
$$A_\tau\stackrel{p}{\to} A_H=A_{\tau/n}\stackrel{q}{\to}A_\tau$$
where $p(z+\Lambda_\tau)=z+\Lambda_{\tau/n}$ and $q(z+\Lambda_{\tau/n})=nz+\Lambda_\tau$. We see that the cosets of $H$ are precisely $\mu+H$ for $\mu\in\frac{1}{n}\mathbb{Z}^g/\mathbb{Z}^g$, and moreover
$$q^{-1}(\mu+\Lambda_\tau)=\frac{1}{n}\mu+\frac{1}{n}\mathbb{Z}^g+\Lambda_{\tau/n}.$$
Then 
$$\varphi_H(q^{-1}(\mu+\Lambda_{\tau}))=\left\{\left[\theta\left[\begin{matrix}\delta\\0\end{matrix}\right](\tau,\mu+a)\right]_{\delta\in\frac{1}{n}\mathbb{Z}^g/\mathbb{Z}^g}:a\in\mathbb{Z}^g/n\mathbb{Z}^g\right\}.$$
But $\theta\left[\begin{matrix}\delta\\0\end{matrix}\right](\tau,\mu+a)=\mbox{exp}(2\pi i \delta^ta)\theta\left[\begin{matrix}\delta\\\mu\end{matrix}\right](\tau,0)$. Therefore,
$$Q_{H,\mu}+1=\mbox{rank}\left(\mbox{exp}(2\pi i n \delta^t\epsilon)\theta\left[\begin{matrix}\delta\\\mu\end{matrix}\right](\tau,0)\right)_{\delta,\epsilon\in\frac{1}{n}\mathbb{Z}^g/\mathbb{Z}^g},$$
and so we have
$$n^{2g}-n^g-Q_H=n^{2g}-\sum_{\mu\in\frac{1}{n}\mathbb{Z}^g/\mathbb{Z}^g}\mbox{rank}\left(\mbox{exp}(2\pi i n \delta^t\epsilon)\theta\left[\begin{matrix}\delta\\\mu\end{matrix}\right](\tau,0)\right)_{\delta,\epsilon\in\frac{1}{n}\mathbb{Z}^g/\mathbb{Z}^g}.$$
A quick check shows that the sum above is equal to the number of non-vanishing theta constants, which we know is equal to $n^{2g}-\Theta(n)$.
\end{proof}

We can now obtain an explicit bound for the number of torsion points on a theta divisor. 

\begin{theorem}
Let $(A,\Theta)$ be a complex principally polarized abelian variety. Then
$$\Theta(2)\leq 4^{g}-g2^{g-1}-2^g$$
and for $n\geq3$
$$\Theta(n)\leq n^{2g}-(g+1)n^g.$$
\end{theorem}
\begin{proof}
By the previous proposition, we can take $\Theta=\Theta_\tau$ on $A_\tau$ for some $\tau\in\mathcal{H}_g$. Using the notation as in the proof of the previous proposition, we have that
$$q^{-1}(\mu+\Lambda_\tau)=\frac{1}{n}\mu+\frac{1}{n}\mathbb{Z}^g+\Lambda_{\tau/n}$$
for $\mu\in\frac{1}{n}\mathbb{Z}^g/\mathbb{Z}^g$. Therefore each member of $\varphi_H(q^{-1}(\mu+\Lambda_\tau))$ differs from the other by the action of the representation $\rho_H:\ker q\to PGL(H^0(A_H,q^*\mathcal{L}))$. It is known that this action (for this particular subgroup) multiplies the projective coordinates of $\mathbb{P}H^0(A_H,q^*\mathcal{L})$ by $n$th roots of unity, and so we can estimate $Q_{H,\mu}$ by the number of vanishing coordinates. Using this fact, it is easy to see that $Q_{H,c_i}$ is equal to $n^g-1-r$, where $r$ is the number of vanishing coordinates. Therefore $n^{Q_{H,\mu}}\geq\#\varphi_H(q^{-1}(\mu))$.

For $n=2$, when $\mu=0$ we have $2^g$ different points in $\varphi_H(\ker q)$, and so $Q_{H,0}\geq g$.

Let us assume $\Theta_{\tau/2}$ irreducible. When  $\mu\neq0$, we have $2^{g-1}$ different points, and so $Q_{H,\mu}\geq g-1$. Adding everything up we get
$$\Theta(2)\leq 4^g-2^g-g-(g-1)(2^g-1)=4^g-g2^g-1.$$
In the case  $\Theta_{\tau/2}$ reducible we proceed in the same  way, but   now we have less points since the map is not injective on the Kummer  variety. The worst case will be when  $(X, \Theta)$ is  a product of elliptic  curves. In this case depending on $\mu$ we can get in the  image $2^{k}$ different points, $k=0, \dots g-1$.  Varying $\mu$ this happens exactly $g\choose k$ times. Hence  totally we  get
$$\Theta(2)\leq 4^g-\sum_{k=0}^g  {g\choose k}(k+1)=4^g-g2^{g-1}-2^g.$$
For $n\geq3$, we have that $\varphi_H$ is an embedding, and so there are always $n^g$ points in $\varphi_H(q^{-1}(c_i))$. This means that $Q_{H,c_i}\geq g$. Therefore if $n\geq3$,
$$\Theta(n)\leq n^{2g}-n^g-gn^g=n^{2g}-(g+1)n^g.$$
\end{proof}

When $\Theta$ is  reducible, even more can be said:

\begin{corollary}
If $(A,\Theta)\simeq\prod_{i=1}^s(B_i,\Theta_i)$ and $b_i=\dim B_i$, then 
$$\Theta(2)\leq 4^{g}-2^g\prod_{i=1}^s\left(\frac{b_i}{2}+1\right)$$
and for $n\geq 3$
$$\Theta(n)\leq n^{2g}-n^{g}\prod_{i=1}^s(b_i+1).$$
\end{corollary}
\begin{proof}
In this case, we see that the number of $n$-torsion points on $\Theta$ is equal to $n^{2g}-t$ where $t$ is the number of $n$-torsion points of the form $(x_1,\ldots,x_s)$ such that $x_i\notin \Theta_i$ for all $i$. Therefore
$$\Theta(2)=4^{g}-\prod_{i=1}^s(4^{b_i}-\Theta_i(2))\leq 4^{g}-\prod_{i=1}^s(b_i2^{b_i-1}+2^{b_i}).$$
The same technique can be applied for $n\geq3$.
\end{proof}

\begin{remark}
If   $(X, \Theta)$ is simple (or more generally not 2-isogenous to a product), using the action  of the symplectic  group
we can improve the  estimate for    $Q_{H,0}$; in fact we   can  get  $Q_{H,0}\geq 2g-1$. Thus in this case we get
$$\Theta(2)\leq 2^{2g}-2^g-g2^g=2^{2g}-(g+1)2^g.$$
 This  number  fits in the general  picture.

\end{remark}

\section{Alternative approaches for $n=2$}

\subsection{Alternative approach 1} 
The methodology in this section is different from that in the previous one, and there are changes in notation. Assume that $\Theta$ is symmetric and irreducible, and define  
$$B_n:=H^0(A,\mathcal{O}_A(n\Theta)).$$
Let $B_n^+$ be the eigenspace associated to 1 for the automorphism $(-1)^*$. It is well-known that 
$$\dim_\mathbb{C}B_n^\pm=2^{g-1}(m^g\pm1).$$
We will use a few results from \cite{K}. For $n\geq 2$ and $m\geq3$, the natural map
$$B_n\otimes B_m\to B_{m+n}$$
is surjective. Since $B_2=B_2^+$, we have that $B_2\otimes B_m^{\pm}\to B_{m+2}^{\pm}$ is surjective, and therefore
$$\mbox{Sym}^2(B_2)\otimes B_m^{\pm}\to B_{m+4}^\pm$$
is surjective. Let $V_2\subseteq B_4^+$ be the image of $\mbox{Sym}^2(B_2)$ in $B_4^+$.  We are interested in a basis  of $V_2$, which is given by all $\theta\left[\begin{matrix}\delta\\ \epsilon\end{matrix}\right](\tau,2z)$ for $\delta,\epsilon\in\frac{1}{2}\mathbb{Z}^g/\mathbb{Z}^g$  and  $4\delta^t\epsilon\equiv 0\mbox{ (mod }2)$ such that $\theta\left[\begin{matrix}\delta\\ \epsilon\end{matrix}\right](\tau,0)\neq0$ (in this section all theta characteristics will be half-integer characteristics). Let $\frak{n}_g$ be the dimension of $V_2$. It is clear that 
$$\Theta(2)=4^g-\frak{n}_g,$$
 since it is the  number of vanishing  theta constants. As an immediate consequence of the previous discussion we have
\begin{proposition}  $$\Theta(2)\leq 4^g-\frac{7^g-1}{3^g-1}$$
\end{proposition}

\begin{proof} We have that the map $\mbox{Sym}^2(B_2)\otimes B_m^{\pm}\to B_{m+4}^\pm$ factors as

\centerline{\xymatrix{\mbox{Sym}^2(B_2)\otimes B_m^{\pm}\ar[dr]\ar[rr]&&B_{m+4}^\pm\\&V_2\otimes B_m^\pm\ar[ur]&}}

\noindent and since the above arrow is surjective, all the arrows are surjective. Therefore,
$$\frak{n}_g\geq \dim_\mathbb{C}B_{m+4}^\pm/\dim_{\mathbb{C}}B_m^\pm=\frac{(m+4)^g\pm1}{m^g\pm1}$$
for $m\geq 3$. The maximum of this function in $m$ is achieved when $m=3$ and the sign is negative.

\end{proof}

\subsection{Alternative approach 2} From the addition  formula for theta functions with semi-integral characteristics  (see \cite[Theorem 2, pg. 139]{Igusa} we have
$$\theta\left[\begin{matrix}\delta\\ \epsilon\end{matrix}\right](\tau,0)\theta\left[\begin{matrix}\delta\\ \epsilon\end{matrix}\right](\tau,2z)=\sum_{\sigma} (-1)^{ <2\epsilon, 2\sigma>}\theta\left[\begin{matrix}\sigma\\0\end{matrix}\right](2\tau,2z) \theta\left[\begin{matrix}\delta+\sigma\\ 0\end{matrix}\right](2\tau,2z).$$

Moreover we can restate this saying that  
$$\Theta(2)-2^{g-1}(2^g-1)=2^{g-1}(2^g+1)-\frak{n}_g$$ 
is the dimension of the space of quadrics that vanish on the image of the \emph{Kummer variety} $K(A)=A/\pm1$, via the embedding $K(A)\hookrightarrow|2\Theta|\simeq\mathbb{P}^{2^g-1}$.

 Since the  Kummer variety is  irreducible and the map is  finite, we have that the image of $K(A)$ cannot be contained in any quadric of  rank 2 in $\mathbb{P}^{2^g-1}$. These quadrics form a  variety of dimension $2^{g+1}-1$ in the space of all quadrics in $\mathbb{P}^{2^g-1}$. Thus we have as a rough estimate:

\begin{lemma} $\frak{n}_g\geq 2^{g+1}-1$.
 \end{lemma}
 
\begin{proof}
The space of  quadrics containing the image of $K(A)$  does not intersect the above described  variety. \end{proof}

This then gives us the bound
$$\Theta(2)\leq 4^g-2^{g+1}+1.$$
This estimate is  very rough and a   careful analysis  could  produce  better results. For example we know that  if $\Theta$  is   irreducible, the  number of   vanishing  quadrics is 
  equal to  $1, 10$ when  $g=3,4$ respectively, and $\geq 66$ when $g=5$.
  All these   are  triangular   numbers that   could  give the dimension  of  the  space of quadrics of bounded rank.\smallskip
  
\subsection{Alternative approach 3} This method is different than the previous approach but gives us the same estimate. We have a short exact sequence
$$0\to R\to V_2\otimes B_4^+\to B_8^+\to0$$
where $R$ is the  space of relations. Let $W_2\subseteq B_4^+$ be such that $B_4^+=V_2\oplus W_2$; it has as a basis the set of theta functions with even characteristics that correspond to a point of order 2 on $\Theta$. Recall that the Heisenberg group, given as a set
$$H=\mathbb{G}_m\times\mathbb{F}_2^g\times\mbox{Hom}(\mathbb{F}_2,\mathbb{G}_m)^g,$$
is a non-commutative group that is non-canonically isomorphic to the theta group of $2\Theta$. Now $H$ acts on $B_4^+$ and $B_8^+$ and decomposes these spaces with respect to its characters. Moreover, the characters are in one to one correspondence with the points of order 2 on $A$. It is known (see \cite[Section 2.4]{vG2}) that for a character $\chi$
\begin{eqnarray}\nonumber\dim(B_8^+)_\chi&=&\left\{\begin{array}{ll}2^g&\mbox{if }\chi\mbox{ is trivial}\\2^{g-1}&\mbox{if not}\end{array}\right.\\
\nonumber\dim(B_4^+)_\chi&=&\left\{\begin{array}{ll}1&\mbox{if }\chi\mbox{ corresponds to an even characteristic}\\0&\mbox{if not}\end{array}\right..
\end{eqnarray}

\begin{lemma}
We have an exact sequence
$$0\to R_0\to\bigoplus_{\chi}(V_2)_\chi\otimes(V_2)_\chi\to(B_8^+)_0\to0,$$
where the subscript 0 refers to the eigenspace corresponding to the trivial character.
\end{lemma}
\begin{proof}
This follows from the surjectivity of $\mbox{Sym}^2(V_2)\oplus(W_2\otimes V_2)\to B_8^+$.
\end{proof}

\begin{corollary}
We have $\frak{n}_g=2^g+\dim R_0$; or in other words, $\Theta(2)=4^g-2^g-\dim R_0$.
\end{corollary}

In order to estimate $\Theta(2)$, we need a better grasp on what $R_0$  or a suitable subspace is. Denote by $K_g^+$ and $K_g^-$ the sets of isotropic (respectively anisotropic) elements in $\mathbb{F}_2^{2g}$ with respect to the quadratic form
$$\langle X,X\rangle=x_1x_{g+1}+\cdots+x_gx_{2g},$$
and let $k_g^+$ and $k_g^-$ be their orders. We introduce the matrix
$$M(g)=M:=\left(\exp\left[i\pi\sum_{i=1}^g(m_in_{g+i}-n_im_{g+i})\right]\right)_{m,n\in\mathbb{Z}^{2g}/2\mathbb{Z}^{2g}}.$$
Now $M$ has the decomposition
$$M=\left(\begin{array}{cc}M^+&N\\N^t&M^-\end{array}\right)$$
where $M^+$ (respectively $M^-$) is the submatrix of $M$ given by the restriction to $K_g^+\times K_g^+$ (respectively $K_g^-\times K_g^-$). The following proposition is proven in \cite[Lemma 1.1]{Fay}:

\begin{proposition}$M$ has two eigenspaces of  dimension  $k_g^+$
and
$k_g^-$ with eigenvalues $\pm 2^g$, while $M^{\pm}$ has eigenspaces of
dimension $(1/3) (2^g \pm 1)(2^{g-1} \pm 1)$ and $(1/3)(2^{2g}-1)$ with
eigenvalues $\pm 2^g$ and $\mp 2^{g-1}$. For $X\in \Bbb C^{k_g ^+}$ and
$Y\in
\Bbb C^{k_g ^-}$, we have
\begin{eqnarray}\nonumber M\begin{pmatrix} X\\ Y\end{pmatrix}=  2^g\begin{pmatrix} X\\ Y\end{pmatrix} &\iff& M^-
Y= 2^{g-1}Y=N^tX\\
\nonumber M\begin{pmatrix} X\\ Y\end{pmatrix}=  -2^g\begin{pmatrix} X\\ Y\end{pmatrix}  &\iff& M^+
X=-2^{g-1}X=NY\\
\nonumber M^+ X= 2^g X &\iff& N^tX=0\\ 
\nonumber M^- Y=-2^gY&\iff& NY=0\\
\nonumber M^+ X=-2^{g-1} X &\mbox{if}& M^+ X-NY=0\\ 
\nonumber M^- Y=2^{g-1}Y &\mbox{ if }& N^tX-M^-Y=0
\end{eqnarray}
\end{proposition}
 If $m=(a,b)\in K_g^+$ for $a$ and $b$ considered as elements of $\{0,1\}^g$, then we use the notation $\theta_m(\tau,z):=\theta\left[\begin{matrix}a/2\\b/2\end{matrix}\right](\tau,z)$. The following lemma is also proved in \cite{Fay}:
 
\begin{lemma} If $X=(v_m)_{m\in K_g^+}$ is a column of $N$, then $M^+ X=-2^{g-1} X$.
 Moreover we have
$$\sum_{m\in K_g^+}v_m\theta_m(\tau,0)^2\theta_m(\tau,2z)^2=0$$
where $(v_m)_{m\in K_g^+}$ is a column of $N$. \end{lemma}
Since  we have
 $$B:=NN^t= 2^{g-1} ( 2^g  I- M^+ ),$$
it is easy to deduce that  $\mbox{rk}(N)=\frac{1}{3}(4^g-1)$.  Thus the columns of $N$
span the whole eigenspace of  $M^+$ with  eigenvalue   $-2^{g-1 }.$ If $S_0\subset  R_0$  is  the subspace  spanned  by  these relations, then we have
$$\dim S_0\leq \frac{1}{3}(4^g-1).$$
Obviously the dimension of $S_0$ is $\frac{1}{3}(4^g-1)$ if there are  no theta constants vanishing. If there are theta constants that vanish then the dimension could drop.

Let $J$ be the $k_g^+\times k_g^+$  diagonal matrix   whose entries are $0$ or $1$ depending on whether or not the theta constant $\theta_m(\tau,0)$   corresponding to $m\in K_g^+$ vanishes. We see that
$$\dim S_0=\mbox{rk}(JN)=\mbox{rk}(JN(JN)^t)=\mbox{rk}(JBJ^t)$$
where $B=NN^t$.  Now  deleting the $0$ rows and columns,  $JBJ^t$ corresponds to a certain principal submatrix  $B_r$ of $B$ of size $r\times r$ where
$$r=\frak{n}_g\geq 2^g+\dim S_0=2^g+\mbox{rk}(JBJ^t).$$
Thus  to  have an estimate for $\frak{n}_g$, we need to estimate the ranks of principal submatrices of $B$. We therefore obtain:

\begin{proposition}
$$\Theta(2)\leq 4^g-2^g-h_0$$
where $h_0=\min\{k \geq 2^g+\mbox{rank}(S):S\mbox{ principal submatrix of }B\mbox{ of order k }   \}$.
\end{proposition}

\begin{corollary}\label{minors} $$\Theta(2)\leq 4^g-2^{g+1}+1$$
 \end{corollary}
\begin{proof}
We will show that all principal minors of $B$ of size $s\leq 2^g-1$ are positive definite. The matrix $B_{2^g -1} $ is semi-positive definite. The entries  are equal to  $2^g-1$ along the diagonal and $\pm 1$ out of the diagonal.   For every $X\in \mathbb{R}^{2^g-1}$ we have
$$  X^tB_{2^g -1} X=\sum_{ 1\leq i<j\leq k} (x_i\pm x_j)^2 + \sum_{ i=1}^{2^g-1} x_i^2.$$
Thus it is positive definite.
 \end{proof}
 
Now $\mbox{Sp}(2g,\mathbb{F}_2)$ acts on the set of characteristics by
$$\left(\begin{array}{cc}a&b\\c&d\end{array}\right)\cdot\left[\begin{matrix}\delta\\ \epsilon\end{matrix}\right]:=\left(\begin{array}{cc}d&-c\\-b&a\end{array}\right)\left(\begin{matrix}\delta\\ \epsilon\end{matrix}\right)+\left(\begin{matrix}\mbox{diag}(c^td)\\\mbox{diag}(a^tb)\end{matrix}\right).$$
This action is double transitive on the set of even (respectively odd) characteristics. Therefore if we want to compute the rank of submatrices of the matrix $B$, we can consider only orbits with respect to the action of this group.

The Kronecker product of $g$ times the matrix $M^+(1)$ is a matrix $L(g)$ of order $3^g$ with eigenvalues $(-1)^k2^{g-k}$ that have multiplicity $\binom{g}{k}2^{g-k}$ for $k=0,\ldots,g$. If we look at the submatrix $B_k$ indexed by all even characteristics $m=\left[\begin{matrix}\delta\\ \epsilon\end{matrix}\right]$ satisfying $4\delta^t\epsilon=0$ in $\mathbb{Z}$, then
$$B_k=2^{g-1}(2^gI_{3^g}-L(g))$$
and has rank $3^g-2^g$. We see that this implies the well-known result that if $(A,\Theta)$ is the product of elliptic curves, then there are $3^g$ points of order two that are not on $\Theta$.

We finish our analysis by looking at the genus 2 case. Double transitivity of the action of the symplectic group implies that
all submatrices of degree 8 of $M^+(2)$ are conjugate via the action
of the symplectic group. For one of these matrices, we can prove that the rank is 5. This implies that 
$$\frak{n}_2\geq 9,$$
which is sharp. We therefore conjecture the following that would imply that $\Theta(2)\leq 4^g-3^g$ for all $g$:

\begin{conj}
The number $h_0$ is reached at $L(g)$.
\end{conj}

\end{document}

\maketitle
\begin{abstract}
In this  note we revisit  a method used in \cite{APS} to  give a sharp  bound  for the number of 2-torsion points on a theta divisor.
\end{abstract}

\section{Introduction}

Let $A$ be a complex abelian variety and let $\mathcal{L}=\mathcal{O}_A(\Theta)$ be a principal polarization on $A$. We set  
$$\Theta(n):=\#A[n]\cap \Theta,$$
where $A[n]$ is the group of $2$-torsion points on $A$. 

In \cite{PM}, the  authors  gave   a bound   for the number of 2-torsion points on a theta divisor. This  bound  has been recently improved in  \cite{APS} where also a bound for the $n$  torsion points  is given.  However, these bounds were  not  optimal . In   \cite{PM}, it  has been    conjectured that the  bound for the  two  torsion  points  is $4^g-3^g$ and is achieved if and only if $(A,\mathcal{L})$ is the polarized product of elliptic curves.  Similarly  in  \cite{APS}  this  conjecture has been   generalized  to the case of  $n$-torsion points.  In these cases the  bound is  $n^{2g} -(n^2-1)^g$. The aim of this  note is  to   prove the first  part   of the first conjecture and we will give an estimate for  $\Theta(n)$  when  $n$ is even , more exacly we  will  prove 
the following:

\begin{theorem} \label{main}
Let $(A,\Theta)$ be a principally polarized abelian variety. Then
$$\Theta(2)\leq 4^{g}- 3^g$$

$$\Theta(2m)\leq m^{2g}(4^g-3^g)$$

\end{theorem}

 The  proofs are    consequence of results proved       in \cite{K1} .

\section{The proof}
We recall shortly  some  basic facts. Let $\tau\in\mathcal{H}_g$ be a matrix in the Siegel upper-half space, and for $\delta,\epsilon\in\mathbb{R}^g$ and $z\in\mathbb{C}^g$ define the theta function with characteristics
$$\theta\left[\begin{matrix}\delta\\ \epsilon\end{matrix}\right](\tau,z):=\sum_{m\in\mathbb{Z}^g}\mbox{exp}[\pi i(m+\delta)^t\tau (m+\delta)+2\pi i(m+\delta)^t(z+\epsilon)].$$
When $\delta=\epsilon=0$ we obtain \emph{Riemann's theta function} $\theta(\tau,z):=\theta\left[\begin{matrix}0\\ 0\end{matrix}\right](\tau,z)$; the projection of $\{\theta(\tau,\cdot)=0\}$ to $A_\tau:=\mathbb{C}^g/\tau\mathbb{Z}^g+\mathbb{Z}^g$ gives  the symmetric theta divisor $\Theta_\tau$. We set $\mathcal{L}_\tau:=\mathcal{O}_{A_\tau}(\Theta_\tau)$. A basis for  
$H^0(A_\tau,\mathcal{L}_\tau^{n^2})$ is  given by 
 $$\left\{\theta\left[\begin{matrix}\delta\\ \epsilon\end{matrix}\right](\tau,nz):\delta, \epsilon\in\frac{1}{n}\mathbb{Z}^g/\mathbb{Z}^g\right\}.$$
 We shortly recall some  basic  facts about these  functions. The evaluation at $z=0$ of the above   functions vanishes if and only if  if the  $n$  torsion   $\tau\delta+ \epsilon$  belongs to $\Theta_{\tau}$. 
 
 It is a well  known  fact  that the dual torus $\hat{A}$ parametrizes isomorphism classes of line  bundles in  $Pic^0(A)$.   The principal polarization  $\Theta_{\tau}$  induces  an   isomorphism  $A\to\hat{A}$  sendindg a point $x$ to the line  bundle $L_x=t^*\mathcal L\otimes\mathcal L^{-1} $ in $\hat{A}$.
 We recall  from \cite{Kempf} the following fundamental  result  (Theorem 6.8), which also appears in \cite[Proposition 1.5]{S} and \cite[Proposition 7.2.2]{BL}.

\begin{theorem} Let $\mathcal M$ be an  ample line  bundle  on $A$
\begin{itemize}
\item For  $L$ in a non-empty subset of $\hat{A}$ the    multiplication
$$H^0(A, \mathcal M^{\otimes 2}\otimes L\otimes M)\otimes H^0(A, \mathcal  M^{\otimes 2}\otimes N\otimes M^{-1})\to  H^0(A, \mathcal M^{\otimes 4}\otimes L\otimes N )$$
 is  surjective  for  fixed $M$  and $N$  in $\hat{A}$.

\item For  $L$ in a non-empty subset of $\hat{A}$ the    multiplication
$$H^0(A, \mathcal M^{\otimes 2}\otimes M )\otimes H^0(A, \mathcal  M^{\otimes 2}\otimes N\otimes L)\to  H^0(A, \mathcal M^{\otimes 4}\otimes L\otimes M\otimes N )$$
 is  surjective  for  fixed $M$  and $N$  in $\hat{A}$.

\item If $n\geq 2$ and $m\geq 3$ then  the  multiplication
$$H^0(A, \mathcal M^{\otimes m}\otimes M )\otimes H^0(A, \mathcal  M^{\otimes n}\otimes N  )\to  H^0(A, \mathcal M^{\otimes n+m}\otimes M\otimes N )$$
 is  surjective  for  arbitrary $M$  and $N$  in $\hat{A}$.
\end{itemize}
\end{theorem}
 
As an immediate   consequence  we get the  following:
\begin{corollary} \label{cor} Let $\mathcal M$ be an  ample line  bundle  on $A$. For  $L$ in a non-empty subset of $\hat{A}$ the    multiplication
$$H^0(A, \mathcal M^{\otimes 2}  )\otimes H^0(A, \mathcal M^{\otimes 2}  )\otimes H^0(A, \mathcal  M^{\otimes 2}\otimes L)\to  H^0(A, \mathcal M^{\otimes 6}\otimes L )$$
is  surjective.
\end{corollary}
\textit{Proof:} It is enough in the second item of the  previous theorem to  take $M=N= \mathcal O_A$ and then apply the third item. \bigskip

Relatively to the  map 

$$ H^0(A, \mathcal M^{\otimes 2}  )\otimes H^0(A, (\mathcal M\otimes L_{2y})^{\otimes 2})\to  H^0(A, (\mathcal M\otimes L_y)^{\otimes 4})$$
 we know  from \cite{K1} that for $\mathcal M=\mathcal{L}_\tau$ the    dimension of the image,  say $m(0,2y)$ is equal to the   number of  $\eta \in A[2]$
such that $2y+\eta\notin\Theta$. This is equivalent to the number of  non-vanishing
  $$\left\{\theta\left[\begin{matrix}\delta\\ \epsilon\end{matrix}\right](\tau,2y):\delta, \epsilon\in\frac{1}{2}\mathbb{Z}^g/\mathbb{Z}^g\right\}.$$

 In  particular  we  have that the  map is  surjective  once  $$2y\notin \bigcup_{\eta \in A[2]} (\Theta+\eta), $$ 
i.e.  no theta function with half integral characteristic  vanishes at the point $2y$. Moreover   we  have 

\begin{proposition}  For any $y\in A$,  we have    
$$ 3^g\leq m(0,2y)\leq 4^g.$$

\end{proposition}

\textit{Proof:}   The upper bound is  obvious. The lower  bound is an immediate  consequence of  Corollary \ref{cor}. \bigskip

 \noindent Now    we can  prove the  theorem stated in the introduction. \smallskip

\noindent  {\it Proof of Th.\ref{main}.}
Assuming   $y=0$ we get that $$\Theta(2)= 4^g- m(0,0)\leq  4^g-3^g$$

Now  let $n=2m$ even, obviously      $A[2m]/A[2]\equiv (\bZ/m\bZ)^{2g}$.  Now   for each $2y$ in  the   above class the  previous estimate holds. Taking the union  on all representatives we  get
$$ \Theta (n)\leq  n^{2g}- m^{2g} 3^g= m^{2g}(4^g- 3^g)$$
and the theorem is  proved.